\documentclass[11pt]{amsart} 

\usepackage{amsrefs}

\usepackage{geometry} 
\usepackage{booktabs} 
\usepackage{array} 
\usepackage{paralist} 
\usepackage{verbatim} 
\usepackage{subfig} 
\usepackage{mathtools}
\usepackage{amsfonts}
 \usepackage{amsthm}

\theoremstyle{plain}

\newtheorem*{theorem}{Theorem}
\newtheorem*{proposition}{Proposition}
\newtheorem*{question}{Question}

\theoremstyle{definition}
\newtheorem*{definition}{Definition}
\usepackage{fancyhdr} 
\newcommand{\field}{k}

\title{The Nodal Cubic is a Quantum Homogeneous Space}
\author{Ulrich Kr\"ahmer and Angela Tabiri}
\address{University of Glasgow, School of Mathematics
and Statistics, 15 University Gardens, G12 8QW Glasgow,
UK}
\email{ulrich.kraehmer@glasgow.ac.uk,a.tabiri1@research.gla.ac.uk}
\date{} 

\begin{document}
\begin{abstract}
The cusp was recently shown to admit the structure of a
quantum homogeneous space, that is, its coordinate ring
$B$ can be embedded as a right coideal subalgebra into
a Hopf algebra $A$ such that $A$ is faithfully flat as
a $B$-module. In the present article such a Hopf
algebra $A$ is constructed for the coordinate ring $B$
of the nodal cubic, thus further motivating the question which
affine varieties are quantum homogeneous spaces.
\end{abstract}
\maketitle

Just as quantum groups (Hopf algebras) generalise
affine algebraic groups, quantum homogeneous spaces
as studied e.g.~in
\cite{MR3473433,MR1992890,MR3144627,MR2944442,MR2931396,MR1710737,MR1098988,MR2872439,MR2802552,MR3428362}
generalise affine varieties with a transitive action of
an algebraic group: 

\begin{definition}
A \emph{quantum homogeneous space} is a 
right coideal subalgebra $B$ of a Hopf algebra $A$
which is faithfully flat as a left $B$-module.
\end{definition}   
\vspace{-\parskip}

There is also an analytic theory of transitive or
more generally ergodic actions of compact or
locally compact quantum groups, see
e.g.~\cite{MR3121622,MR3128415} and
the references therein.

The best studied examples
are deformation
quantisations of affine homogeneous spaces, and in
particular Podle\'s' quantum 2-spheres \cite{MR919322}. 
However, even if $A$ is noncommutative, $B$ can be a
commutative algebra, so a natural question to ask is:

\begin{question}
Which  
affine varieties are quantum homogeneous spaces?
\end{question}

\vspace{-\parskip}

In the operator algebraic setting, the analogue of this
question has been raised and studied for example in
\cite{MR3054297}. Here we consider the purely algebraic
setting, working over an algebraically
closed field $k$ of characteristic 0. 

Unlike a 
homogeneous space, an affine variety which is a quantum
homogeneous space
can be singular as the example of 
the cusp shows (see \cite[Section~2.11]{MR2931396} and
\cite[Construction~1.2]{MR2732991}). Note that by
the results from \cite{etingof}, noncommutative Hopf
algebras coacting on commutative algebras are quite
restricted. Still, we conjecture that every plane curve 
can be given the structure of a quantum homogeneous
space, and our aim here is to
point this out for the nodal cubic given by the equation
$y^2=x^2+x^3$:

\begin{theorem}\label{accra}
Fix
$(q,p) \in \field^2$ satisfying $p^2=q^2+q^3$. Then the
unital associative $\field$-algebra $A$ with generators 
$x,y,a,a^{-1},b$
satisfying the relations
\begin{gather*}
	a a^{-1} = a^{-1} a = 1,\qquad
	y^{2} = x^{2} + x^{3},\qquad
	b^{2} 
	= a^{3},\\
	ba = ab,\qquad 
	ya = ay,\qquad 
	bx = xb,\qquad
	yx = xy,\qquad
	by = -yb+2pb^{2},\\
	a^2 x 
= - xa^{2} - axa - a^2 +
	\left( 1 + 3q  \right) a^{3},\\
	ax^2
= -ax-xa-x^2a-xax+(2+3q)qa^3.
\end{gather*}
admits a Hopf algebra structure whose coproduct
$ \Delta $, counit $ \varepsilon $ and antipode $S$
satisfy
\begin{gather*}
	\Delta(x) 
= 1 \otimes \left( x - qa  \right) + x \otimes
a,\qquad
	\Delta(y) 
= 1 \otimes \left( y - pb  \right) + y \otimes b,\\
	\Delta (a)=a \otimes a,\quad
	\Delta(b) 
= b \otimes b,\quad
	\varepsilon(x) = q, \quad
	\varepsilon(y)=p, \quad
	\varepsilon(a)= \varepsilon(b)=1,\\
	S(x) = q - \left( x - q \right)a^{-1},\quad
	S(y) = p - \left( y - p \right)b^{-1},\quad
	S(a) = a^{-1},\quad 
	S(b) = b^{-1}.
\end{gather*}
Furthermore, the right coideal subalgebra $B \subset A$ generated by $x,y$ is
the coordinate ring of the nodal cubic, and $A$ is free
and in particular faithfully flat as a $B$-module. 
\end{theorem} 
\vspace{-\parskip}

Observe that the commutation relations in $A$ are chosen  
in such a way that 
$$
	(y-pb)^2=y^2-p^2b^2,\qquad
	(x-qa)^2+(x-qa)^3=
	x^2+x^3-(q^2+q^3)a^3
$$
so that
$$
	(y-pb)^2=(x-qa)^2+(x-qa)^3.
$$
Thus informally speaking each coordinate on the curve 
becomes perturbed by some additional group-like ``quantum
coordinate'' and the perturbed coordinates still
satisfy the defining relation of the curve. 
As $x-qa$ and $y-pb$ are twisted
primitive, the Hopf algebra $A$ is generated by group-likes
and twisted primitives. This implies:

\begin{proposition}
The Hopf algebra $A$ is pointed.
\end{proposition}
\vspace{-\parskip}

The point 
$(q,p)$ on the curve is the point the quantum orbit of 
which it is presented as. To use these observations as starting
point of a general study of quantum homogeneous space
structures on affine varieties seems 
a promising future research direction.  

The proof of the theorem consists of a straightforward
(albeit tedious) verification that the formulas
for the coproduct, counit and antipode are compatible 
with the defining relations of $A$, followed by a
similarly straigthforward application of Bergman's diamond
lemma \cite{MR506890}
yielding a vector space basis of $A$ that implies the
freeness over $B$:

\begin{proposition}
The set
$$
	\{ x^i y^j (ax)^l a^m b^n \mid 
	i,l \in \mathbb{N} , j \in \{0,1\}, 
	 m \in \mathbb{Z}, n \in \{0,1\}\}
$$
is a vector space basis of $A$,
and the GK-dimension of $A$ equals 3.
\end{proposition}
\vspace{-\parskip}

Using this basis, one also easily observes 
that like the nonstandard Podle\'s spheres, 
the algebra extension $B \subset A$ is
an example of a coalgebra Galois extension
\cite{MR1684158,MR1669095,MR2038278,MR1604340} 
rather than of a Hopf-Galois extension: 
the
coalgebra is $C:=A/B^+A$, where $B^+:=B \cap
\mathrm{ker}\,\varepsilon $.
The canonical projection $ \pi : A \rightarrow C$
defines a left $C$-coaction
$$
   \lambda \colon A \rightarrow C \otimes A,\quad
   f \mapsto f_{(-1)} \otimes f_{(0)} := 
   \pi (f_{(1)}) \otimes f_{(2)} 
$$
and we have (as a consequence of the faithful flatness
of $A$ over $B$)
$$
   B = \{ f \in A \mid 
   f_{(-1)} \otimes f_{(0)} = 
   \pi (1) \otimes f\}.
$$
That $C$ is not a Hopf algebra
quotient of $A$ follows from  
$B^+A \neq A B^+$
(cf.~\cite[Lemma~1.4]{MR1710737}); for example,
we have
$$
   AB^+ \ni a^2(x-q)=
   -xa^2-axa-(1+q)a^2+(1+3q)a^3 \notin B^+A. 
$$

We finally remark that 
some properties of the algebra $A$ are better understood 
when using a slightly different set of generators: if
we abbreviate 
$$
	c:=3 x -(1+3q) a + 1,\quad
	d:=3y-6pb,\quad
	e:=ac+rca
$$
where $r$ is a primitive 6th root of 1 (so that 
$r+r^{-1}=1$), 
then the defining relations of $A$ in terms of the
generators $a^{\pm 1},b,c,d,e$ read
$$
	aa^{-1}=a^{-1}a=1,\quad
	ab=ba,\quad
	ac+rca=e,\quad
	ad=da,\quad
	ae+r^{-1}ea=0,
$$
$$
	bc=cb,\quad
	bd=-db,\quad
	be=eb,\quad b^2=a^3,\quad
	cd=dc,\quad
	r^{-1}ce+ec=3(a-a^3),
$$
$$
	de=ed,\quad 3d^2=c^3-3c+2+(1+3q)(-2+6q+9q^2)a^3.
$$

Using these generators, one easily verifies for example:

\begin{proposition}\label{properties}
The units in $A$ are of the form
$\alpha a^ib^j$, $ \alpha \in k$,
$i \in \mathbb{Z} $, $j \in \{0,1\}$.   
\end{proposition}

\subsection*{Acknowledgements}
We thank Ken Brown and Tomasz Brzezi\'nski for helpful
comments on a first draft of this paper.  
A.T.~is funded by a Faculty for the Future Fellowship of
the Schlumberger Foundation.

\end{document}